\documentclass[12pt,reqno]{amsart}
\usepackage{ amsmath , amssymb }
\usepackage[mathscr]{eucal}

\newdimen\unit\newdimen\psep\newcount\nd\newcount\ndx\newbox\dotb\newbox\ptbox
\newdimen\dx\newdimen\dy\newdimen\dxx\newdimen\dyy\newdimen\hgt
\newdimen\xoff\newdimen\yoff
\newcommand\clap[1]{\hbox to 0pt{\hss{#1}\hss}}
\newcommand\vdisk[1]{{\font\dotf=cmr10 scaled #1\dotf.}}
\newcommand\varline[2]{\setbox\dotb\hbox{\vdisk{#1}}\xoff=-.5\wd\dotb
\wd\dotb=0pt\yoff=-.5\ht\dotb\psep=#2\ht\dotb}
\newcommand\varpt[1]{\setbox\ptbox\clap{\vdisk{#1}}\setbox\ptbox
\hbox{\raise-.5\ht\ptbox\box\ptbox}}
\newcommand\cpt{\copy\ptbox}
\newcommand\point[3]{\rlap{\kern#1\unit\raise#2\unit\hbox{#3}}}
\newcommand\setnd[4]{\dx=#3\unit\advance\dx-#1\unit\divide\dx by\psep
\dy=#4\unit\advance\dy-#2\unit\divide\dy by\psep \multiply\dx
by\dx\multiply\dy by\dy\advance\dx\dy\nd=1\advance\dx-1sp
\loop\ifnum\dx>0\advance\dx-\nd sp\advance\nd1\advance\dx-\nd
sp\repeat}
\newcommand\dl[4]{{\setnd{#1}{#2}{#3}{#4}\dline{#1}{#2}{#3}{#4}\nd}}
\newcommand\dline[5]{{\nd=#5\hgt=#2\unit\dx=#3\unit\advance\dx-#1\unit
\divide\dx by\nd\dy=#4\unit\advance\dy-#2\unit\divide\dy by\nd
\advance\hgt\yoff\rlap{\kern#1\unit\kern\xoff\loop\ifnum\nd>1\advance\nd-1
\advance\hgt\dy\kern\dx\raise\hgt\copy\dotb\repeat}}}
\newcommand\ellipse[4]{\qellip{#1}{#2}{#3}{#4}\qellip{#1}{#2}{#3}{-#4}%
\qellip{#1}{#2}{-#3}{#4}\qellip{#1}{#2}{-#3}{-#4}}
\newcommand\qellip[4]{{\setnd{0}{0}{#3}{#4}\dx=\unit\dy=0pt\raise\yoff\rlap{%
\kern#1\unit\kern\xoff\raise#2\unit\hbox{\loop\ifnum\dx>0\rlap{\kern#3\dx
\raise#4\dy\copy\dotb}\hgt=\dx\divide\hgt
by\nd\advance\dy\hgt\hgt=\dy \divide\hgt
by\nd\advance\dx-\hgt\repeat\rlap{\raise#4\dy\copy\dotb}}}}}
\newcommand\bez[6]{{\setnd{#1}{#2}{#3}{#4}\ndx=\nd\setnd{#3}{#4}{#5}{#6}
\ifnum\ndx>\nd\nd=\ndx\fi\dx=#3\unit\advance\dx-#1\unit\dy=#4\unit
\advance\dy-#2\unit\dxx=#5\unit\advance\dxx-#1\unit\dyy=#6\unit\advance
\dyy-#2\unit\advance\dxx-2\dx\advance\dyy-2\dy\divide\dxx
by\nd\divide\dyy
by\nd\advance\dx.25\dxx\advance\dy.25\dyy\divide\dx
by\nd\divide\dy by\nd \multiply\nd
by2\dx=100\dx\dy=100\dy\dxx=100\dxx\dyy=100\dyy\divide\dxx by\nd
\divide\dyy
by\nd\hgt=#2\unit\raise\yoff\rlap{\kern#1\unit\kern\xoff
\raise\hgt\copy\dotb\loop\ifnum\nd>0\advance\nd-1\advance\hgt0.01\dy
\kern0.01\dx\raise\hgt\copy\dotb\advance\dx\dxx\advance\dy\dyy\repeat}}}
\newcommand\ptu[3]{\point{#1}{#2}{\cpt\raise1ex\clap{$\scriptstyle{#3}$}}}
\newcommand\ptd[3]{\point{#1}{#2}{\cpt\raise-1.8ex\clap{$\scriptstyle{#3}$}}}
\newcommand\ptr[3]{\point{#1}{#2}{\cpt\raise-.4ex\rlap{$\ \scriptstyle{#3}$}}}
\newcommand\ptl[3]{\point{#1}{#2}{\cpt\raise-.4ex\llap{$\scriptstyle{#3}\ $}}}
\newcommand\ptlu[3]{\point{#1}{#2}{\raise.8ex\clap{$\scriptstyle{#3}$}}}
\newcommand\ptld[3]{\point{#1}{#2}{\raise-1.6ex\clap{$\scriptstyle{#3}$}}}
\newcommand\ptlr[3]{\point{#1}{#2}{\raise-.4ex\rlap{$\,\scriptstyle{#3}$}}}
\newcommand\ptll[3]{\point{#1}{#2}{\raise-.4ex\llap{$\scriptstyle{#3}\,$}}}

\newcommand\thkline{\varline{1600}{.3}}
\newcommand\medline{\varline{800}{.5}}
\newcommand\thnline{\varline{400}{.6}}

\varpt{2500}\thnline\unit=2em

\newtheorem{thm}{Theorem}
\newtheorem{conjecture}{Conjecture}
\newtheorem*{mader}{Mader's Theorem}

\newtheorem{prob}{Problem}
\newtheorem{lemma}[thm]{Lemma}
\newtheorem{cor}[thm]{Corollary}

\newtheorem{obs}{Observation}
\theoremstyle{definition}\newtheorem{rmk}{Remark}
\theoremstyle{definition}

\def\N{\mathbb{N}}
\def\R{\mathbb{R}}

\def\eps{\varepsilon}
\def\le{\leqslant}
\def\ge{\geqslant}

\begin{document}
\hspace{0cm}\\[-6ex]

\title[Highly connected monochromatic subgraphs]{Highly connected monochromatic subgraphs of multicoloured graphs}

\author{Henry Liu}
\address{Department of Mathematical Sciences\\ The University of Memphis\\ Memphis, TN 38152} \email{henryliu@memphis.edu}

\author{Robert Morris}
\address{Department of Mathematical Sciences\\ The University of Memphis\\ Memphis, TN 38152} \email{rdmorrs1@memphis.edu}

\author{Noah Prince}
\address{Department of Mathematics\\ University of Illinois\\ 1409 W. Green Street\\ Urbana, IL 61801} \email{nprince@math.uiuc.edu}\thanks{The second and third authors were partially supported during this
research by Van Vleet Memorial Doctoral Fellowships.}

\maketitle

\begin{abstract}
We consider the following question of Bollob{\'a}s: given an $r$-colouring of $E(K_n)$, how large a $k$-connected subgraph can we find using at most $s$ colours? We provide a partial solution to this problem when $s = 1$ (and $n$ is not too small), showing that when $r = 2$ the answer is $n - 2k + 2$, when $r = 3$ the answer is $\lfloor \frac{n-k}{2} \rfloor + 1$ or $\lceil \frac{n-k}{2} \rceil + 1$, and when $r - 1$ is a prime power then the answer lies between $\frac{n}{r-1}-11(k^2-k)r$ and $\frac{n-k+1}{r-1} + r$. The case $s \ge 2$ is considered in a subsequent paper~\cite{HNR2}, where we also discuss some of the more glaring open problems relating to this question.
\end{abstract}

\section{Introduction}\label{intro}

A graph $G$ on $n \ge k+1$ vertices is said to be \textit{$k$-connected} if whenever at most $k-1$ vertices are removed from $G$, the remaining vertices are still connected by edges of $G$. It is easy to see that given any graph $G$, either $G$ or $\overline{G}$ (the complementary graph) is connected. A substantial generalisation of this observation, due to Bollob\'as, asks the following question: When we colour the edges of the complete graph $K_n$ with at most $r$ colours, how large a $k$-connected subgraph are we guaranteed to find using only at most $s$ of the colours? In this paper we shall provide a partial answer to this question in the case $s = 1$, and in a subsequent paper~\cite{HNR2} we shall consider the case $s \ge 2$, and in particular the cases $s = 2$, $2s = r$ and $s = \Theta(\sqrt{r})$, where a jump occurs. The majority of the problem, however, remains wide open.

Bollob\'as and Gy\'arf\'as~\cite{BG} observed the following example in the case $r = 2$ and $s = 1$. First partition the vertices of $K_n$ into five classes, four of order $k-1$ (call these $A_1$, $A_2$, $A_3$ and $A_4$) and the fifth containing the remaining $n - 4k + 4$ vertices (call it $B$). Colour the edges between $A_i$ and $B$ red if $i = 1$ or $2$, and blue if $i = 3$ or $4$, and colour the edges between $A_i$ and $A_j$ red if $\{i,j\} \in \{\{1,2\},\{1,3\},\{2,4\}\}$ and blue otherwise ($i \neq j$). Colour the edges inside the blocks arbitrarily. The construction is pictured below (Figure 1) with only the red edges drawn.

\[ \unit = 1.1cm
\medline \ellipse{0.5}{0}{0.5}{0.5} \ellipse{0.5}{2}{0.5}{0.5}
\ellipse{2.5}{0}{0.5}{0.5} \ellipse{2.5}{2}{0.5}{0.5}
\ellipse{-2}{1}{1}{1} \thkline \dl{-1.2}{0.6}{0.1}{0.1}
\dl{-1.2}{1.4}{0.1}{1.9} \dl{0.5}{0.4}{0.5}{1.6}
\dl{0.9}{0}{2.1}{0} \dl{0.9}{2}{2.1}{2} \ptlu{0.5}{-0.2}{k - 1}
\ptlu{2.5}{-0.2}{k - 1} \ptlu{0.5}{1.8}{k - 1} \ptlu{2.5}{1.8}{k -
1} \ptlu{-2}{0.8}{n - 4(k - 1)} \ptlr{3}{2.4}{A_3}
\ptlr{3}{-0.4}{A_4} \ptll{0}{2.4}{A_1} \ptll{0}{-0.4}{A_2}
\ptll{-3.1}{1}{B} \point{-0.7}{-1.4}{\small Figure 1}
\]\\[-2ex]

How large a $k$-connected subgraph using edges of only one colour does this colouring contain? Suppose such a subgraph $H$ uses more than $n - 2k + 2$ of the vertices, and assume (by the symmetry between the given red and blue edges) that the edges of $H$ are coloured red. $H$ must use some vertex $v$ of $A_3 \cup A_4$; suppose $v \in A_3$. But now if we remove the vertices of $A_1 \cap V(H)$ from $H$ (to get $H'$ say) then $v$ is no longer in the same component of $H'$ as any vertex of $(A_2 \cup A_4 \cup B) \cap V(H')$, and if $n \ge 4k - 3$ then such a vertex must exist.

Since $|A_1 \cap V(H)| \le |A_1| = k-1$, we have shown that if $n \ge 4k - 3$, we cannot guarantee finding a monochromatic $k$-connected subgraph on more than $n - 2k + 2$ vertices. Bollob\'as and Gy\'arf\'as conjectured that this example is extremal, i.e., that if $n \ge 4k - 3$ we \textit{can} guarantee finding a monochromatic $k$-connected subgraph on \textit{at least} $n - 2k + 2$ vertices (note that when $n = 4k-4$ the example above (with $A_1$ and $A_2$ blue cliques, and $A_3$ and $A_4$ red cliques) contains no monochromatic $k$-connected subgraph at all, so the conjecture really is the strongest possible). They also gave a short proof of a somewhat weaker statement~\cite{BG}. Using the ideas from their proof, we are able to show that the conjecture holds when $n \ge 13k - 15$. To state this result we shall need a little notation.

Suppose we are given $n,r,s,k \in \N$, and a function $f: E(K_n) \to [r]$, i.e., an $r$-colouring of the edges of $K_n$. We assume always that $n \ge 2$. Given a subgraph $H$ of $K_n$, write $c_f(H)$ for the order of the image of $E(H)$ under $f$, i.e., $c_f(H) = |f(E(H))|$, the number of different colours with which $f$ colours $H$. Now, define $M(f,n,r,s,k) = \max\{|V(H)| : H \subset K_n$, $c_f(H) \le s\}$, the order of the largest $k$-connected subgraph of $K_n$ using at most $s$ colours from $[r]$. Finally, define $m(n,r,s,k) = \min_f\{M(f,n,r,s,k)\}$. Thus, the question of Bollob\'as asks for the determination of $m(n,r,s,k)$ for all values of the parameters. We shall state all our main results in terms of $m(n,r,s,k)$.

Our first result is the following; it is exactly the conjecture of Bollob\'as and Gy\'arf\'as in the case $n \ge 13k - 15$.

\begin{thm}{\label{21k}}
Let $n,k \in \N$, with $n \ge 13k - 15$. Then $$m(n,2,1,k) \; = \; n - 2k + 2.$$
\end{thm}

Unfortunately our method breaks down when $n$ is much smaller than $13k$, and an analysis of the situation for small values of $k$ suggests that a completely different approach may be necessary in this case.

For $r > 2$ the situation becomes a little more complicated. Many years ago\footnote{\noindent We apologise to those readers who do not consider 25 years to be `many'!}, whilst studying a different problem (on hypergraph covering), F{\"u}redi~\cite{ZF} and Gy{\'a}rf{\'a}s~\cite{AG} showed independently that $\frac{n}{r-1} \le m(n,r,1,1) \le \frac{n}{r-1} + r$ whenever $r - 1$ is a prime power, with equality in the lower bound when $(r-1)^2$ divides $n$. In Section~\ref{sectr1k} we shall give a short proof of this result. It is easy to modify the upper bound construction of F{\"u}redi and Gy{\'a}rf{\'a}s to give $m(n,r,1,k) \le \frac{n-k+1}{r-1} + c_{n,k,r}$ when $r-1$ is a prime power, where $c = c_{n,k,r} \le r$ and $c = 0$ when $(r-1)^2$ divides $(n - r(k-1))$ (see Section~\ref{sectr1k}). The next result shows that this upper bound is essentially best possible for these values of $r$.

\begin{thm}\label{r1k}
Let $n,k,r \in \N$, with $r \ge 3$ and $r - 1$ a prime power. Then,
$$\frac{n}{r-1} \: - \: 11(k^2 - k)r \; \le \; m(n,r,1,k) \; \le \; \frac{n-k+1}{r-1} \: + \: r,$$
and moreover, the lower bound holds for all $3 \le r \in \N$. In particular, if $r$ and $k$ are fixed, then $m(n,r,1,k) = \frac{n}{r-1} + o(n)$.
\end{thm}

Finally, we shall determine the function exactly when $r = 3$.

\begin{thm}\label{31k}
Let $n,k \in \N$, with $n \ge 480k$. Then
$$\frac{n-k+1}{2} \; \le \; m(n,3,1,k) \; \le \; \frac{n-k+1}{2} \: + \: 1.$$
\end{thm}

\noindent Moreover, equality holds in the lower bound of Theorem~\ref{31k} if and only if $n + k \equiv 1 \pmod 4$ (see Corollary~\ref{31kexact}).

The rest of the paper is organised as follows. In Section~\ref{sect21k} we shall prove Theorem~\ref{21k}, and in Section~\ref{sectr1k} we shall prove Theorems~\ref{r1k} and \ref{31k}.

\section{The case $r = 2$}\label{sect21k}

Our first task is to prove Theorem~\ref{21k}. Given any 2-colouring $f$ of $E(K_n)$, we write $R$ for the graph on $V(K_n)$ with edge set $f^{-1}(1)$, and $B$ for the graph with edge set $f^{-1}(2)$, so $E(R) \cup E(B) = E(K_n)$. We shall always refer to the colours as `red' and `blue' respectively.

The set of neighbours of a vertex $x$ in a graph $G$ will be denoted by $\Gamma_G(x)$, or just $\Gamma(x)$ when it is clear to which graph we refer, and similarly the degree of $x$ will be denoted $d_G(x)$, or simply $d(x)$. We shall write $G[A]$ for the subgraph of $G$ induced by a set $A \subset V(G)$, and $G - A$ for the graph $G[V(G) \setminus A]$. If $C, D \subset V(G)$ and $C \cap D = \emptyset$, then $G[C,D]$ will denote the bipartite graph, with parts $C$ and $D$, induced by $G$. For any undefined terms, see~\cite{MGT}.

The following simple lemma appeared in~\cite{BG}. We give the proof for the sake of completeness.

\begin{lemma}\label{degs}
In any $2$-colouring of $E(K_n)$ with $d_R(v) \ge 2k-2$ for every $v \in V(K_n)$, either $R$ is $k$-connected or $B$ contains a $k$-connected subgraph on at least $n - k + 1$ vertices.
\end{lemma}

\begin{proof}
If $R$ is not $k$-connected, then $B$ must contain a complete bipartite graph $H$ on $n-k+1$ vertices. Let the part sizes be $i$ and $j$. If $1 \le i \le k-1$, then $j \ge n-2k+2$, and any vertex $v$ in the $i$-set has $d_R(v) \le 2k-3$, a contradiction. Hence $i \ge k$, and similarly $j \ge k$, so $H$ is $k$-connected.
\end{proof}

We shall also need the following easy lemma, which will be useful throughout the entire paper.

\begin{lemma}\label{intersect}
Let $G$ be a bipartite graph with partite sets $M$ and $N$ such that $d(x) \ge k$ for every $x \in M$, and $|\Gamma(y) \cap \Gamma(z)| \ge k$ for every pair $y,z \in N$. Then $G$ is $k$-connected.
\end{lemma}

\begin{proof}
Let $G$ be such a bipartite graph, and let $C$ be any subset of $V(G)$ of size at most $k-1$. We wish to show that $G' = G - C$ is connected. But this is clear, since any two vertices $x,y \in V(G') \cap N$ have a common neighbour in the graph $G'$ (since $|\Gamma_G(x) \cap \Gamma_G(y)| \ge k$ and $|C| \le k-1$), and any vertex $z \in M$ has a neighbour in $V(G') \cap N$, since $d_G(z) \ge k$. The lemma follows.
\end{proof}

Finally, we make a trivial observation.

\begin{obs}\label{addvtx}
Let $G$ be a graph, and $v \in V(G)$. If $G - v$ is $k$-connected and $d(v) \ge k$, then $G$ is also $k$-connected.
\end{obs}

We are now ready to prove Theorem~\ref{21k}.

\begin{proof}[Proof of Theorem~\ref{21k}]
Let $k \in \N$, with $n \ge 13k - 15$. The upper bound, $m(n,2,1,k) \le n - 2k + 2$, follows from the construction described in Section~\ref{intro}. To prove the matching lower bound, let $f$ be a $2$-colouring of $E(K_n)$. We shall find a monochromatic $k$-connected subgraph of $K_n$ on at least $n - 2k + 2$ vertices.

By Lemma~\ref{degs}, we may assume that there exist vertices $x_1,y_1 \in V = V(K_n)$ with $d_R(x_1) \le 2k - 3$ and $d_B(y_1) \le 2k - 3$, as otherwise the lemma gives us a monochromatic $k$-connected subgraph on at least $n-k+1$ vertices. We construct (by choosing vertices one by one) maximal subsets $X = \{x_1, \ldots, x_p\}$ and $Y = \{y_1, \ldots, y_q\}$ of $V$ such that\\[+1ex] (a) for each $i \in [p]$, $d(x_i) \le 2k-3$ in the graph $R - \{x_1,\ldots,x_{i-1}\}$, and\\[+1ex] (b) for each $i \in [q]$, $d(y_i) \le 2k-3$ in the graph $B - \{y_1, \ldots, y_{i-1}\}$.\\[+2ex]
\underline{Claim 1}: $\min(p,q) \le 8k-11$.

\begin{proof}
Let $u = |X \setminus Y|$, $v = |Y \setminus X|$ and $r = |X \cap Y|$, so $p = u + r$ and $q = v + r$. Let $e_R(X)$ be the number of red edges in $H = K_n[X \cup Y]$ with an endpoint in $X$, and $e_B(Y)$ be the number of blue edges in $H$ with an endpoint in $Y$. Now, there are $uv + ur + vr + {r \choose 2}$ edges in $H$ with an endpoint in $X$ and an endpoint in $Y$. Each such edge contributes at least one to $e_R(X)$ or $e_B(Y)$, so $e_R(X) + e_B(Y) \ge uv + ur + vr + {r \choose 2}$. Also, by the definition of $X$, the number of edges in $R$ with an endpoint in $X$ is at most $(2k - 3)|X|$, so $e_R(X) \le (2k - 3)p$, and similarly $e_B(Y) \le (2k - 3)q$. Hence,
\begin{eqnarray*}
pq & = & uv \: + \: ru \: + \: rv \: + \: r^2 \; \le \; e_R(X)
 \: + \: e_B(Y) \: + \: r^2 \: - \: {r \choose 2}\\
& & \;\;\;\;\;\; \le \; (2k-3)p \: + \: (2k-3)q \: + \:
\frac{1}{2}r^2
\: + \: \frac{1}{2}r\\[+1ex]
& & \;\;\;\;\;\; \le \; (2k-3)(p+q) \: + \: \frac{1}{2}pq \: + \:
\frac{1}{4}(p+q)
\end{eqnarray*}
since $r \le p$ and $r \le q$. It follows that
$$\frac{1}{2}pq \; \le \; \left(2k - \frac{11}{4} \right)(p + q),$$ and so dividing by $pq/2$, we get
$$1 \; \le \; \left( \frac{1}{p} \: + \: \frac{1}{q} \right) \left( 4k - \frac{11}{2} \right).$$
Therefore $p$ or $q$ is at most $8k - 11$.
\end{proof}

Assume then, without loss of generality, that $|X| \le 8k-11$. Note that $X$ was chosen to be maximal, so $d_{R-X}(v) \ge 2k - 2$ for every vertex $v \in V \setminus X$. Therefore, by Lemma~\ref{degs}, either $R - X$ is $k$-connected, or $B - X$ contains a $k$-connected subgraph $H$ on at least $n - |X| - k + 1$ vertices. Suppose the latter. By the definition of $X$, any vertex $x \in X$ sends at most $2k - 3$ red edges into $H$, and so $x$ must send at least
$$|H| - 2k + 3 \; \ge \; n - |X| - 3k + 4
\; \ge \; n - 11k + 15 \; > \; k$$ blue edges into $H$. So by
Observation~\ref{addvtx}, $B[V(H) \cup X]$ is $k$-connected, and
has $n - k + 1$ vertices. Hence we may assume that $R - X$ is
$k$-connected.

Now choose a set $M'$ containing $V \setminus X$ by repeatedly moving from $X$ to $M'$ those vertices which send at least $k$ red edges to $M'$. To be precise, set $X_0 = X$ and $M_0 = V \setminus X$, and at time $t \in \N$ form $X_t$ and $M_t$ by choosing a vertex $v \in X_{t-1}$ with $|\Gamma_R(v) \cap M_{t-1}| \ge k$ if one exists, and setting $X_t = X_{t-1} \setminus \{v\}$ and $M_t = M_{t-1} \cup \{v\}$. If no such vertex exists then stop the process, and set $N = X_{t-1}$ and $M' = M_{t-1}$. Notice that every vertex in $N$ sends at most $k-1$ red edges into $M'$, and that $R[M']$ is $k$-connected by Observation~\ref{addvtx}.

If $|N| \le 2k-2$ then $R[M']$ is our desired subgraph, so assume that $|N| \ge 2k-1$. We wish to apply Lemma~\ref{intersect} to the bipartite graph $G' = B[M',N]$, but we may have some `bad' vertices $v \in M'$ with $d_{G'}(v) \le k-1$. We must therefore first remove these vertices from $M'$.

Let $U$ denote the set of bad vertices in $M'$, so $$U = \{ v \in M' : d_{G'}(v) \le k-1\}.$$ Since each vertex of $N$ sends at most $k-1$ red edges into $M'$, $R[M',N]$ has at most $|N|(k-1)$ edges. But each vertex of $U$ sends at least $|N|-k+1$ red edges into $N$. Thus we have
$$|U|(|N|-k+1) \; \le \; |N|(k-1),$$
and hence
$$|U| \; \le \; \frac{|N|(k-1)}{|N|-k+1} \; \le \; \frac{(2k-1)(k-1)}{k} \;
\le \; 2k-2,$$ since the function $\frac{x(k-1)}{x-k+1}$ is decreasing for $x > k-1$, and $|N| \ge 2k-1$.

We complete the proof of Theorem~\ref{21k} by setting $M = M' \setminus U$, and applying Lemma~\ref{intersect} to the graph $G = B[M,N]$. By the definition of $U$, $d_G(x) \ge k$ for every vertex $x \in M$. Also,
$$|M| \; \ge \; n - |X| - |U|
\; \ge \; n - 10k + 13 \; \ge \; 3k - 2,$$ since $|X| \le 8k -
11$, $|U| \le 2k - 2$ and $n \ge 13k - 15$, and as observed
earlier, $d_G(y) \ge |M| - k + 1$ for every $y \in N$. Therefore,
$$|\Gamma_G(y) \cap \Gamma_G(z)| \; \ge \;
|M| - 2k + 2 \; \ge \; k$$ for every pair $y,z \in N$, so by
Lemma~\ref{intersect}, $G$ is $k$-connected.

Since $M \cup N = V \setminus U$ and $|U| \le 2k-2$, $G$ is the desired monochromatic $k$-connected subgraph.
\end{proof}

\begin{rmk}\label{alpha}
We can in fact improve (for $k \ge 18$) the bound on $n$ to $n \ge (9 + \sqrt{10})k$ as follows. First note that $n - 11k + 15 > k$ still holds, so it will suffice to show that $|M| \ge 3k - 2$. Set $\alpha = 4 + \sqrt{10}$. We have $|M| = n - |N| - |U|$, so if $|N| \le \alpha k + 4$, then $|M| \ge 3k - 2$ if $n \ge (\alpha + 5)k$. If $|N| \ge \alpha k+5$ however, then $$|U| \; < \; \frac{(\alpha k + 5)k}{(\alpha - 1)k + 6} \; < \; \frac{\alpha k}{\alpha - 1} \; = \; (\sqrt{10} - 2)k,$$ so if $n \ge (\alpha + 5)k - 13$ then $|M| \: \ge \: n - 8k + 11 - |U| \: \ge \: 3k - 2$.
\end{rmk}

We have the following rather weak corollary to Theorem~\ref{21k} (and Remark~\ref{alpha}), which would be improved by further reducing the bound on $n$.

\begin{cor}
For every graph $G$ on $n$ vertices, $G$ or $\overline{G}$ has a $\left\lfloor n/(9 + \sqrt{10}) \right\rfloor$-connected subgraph on at least $n - 2 \left\lfloor n/(9 + \sqrt{10}) \right\rfloor + 2$ vertices.
\end{cor}

What happens when $n$ is much smaller? For $n$ close to $4k - 3$ the problem seems to become much more complicated, so we have been forced to restrict ourselves to small values of $k$. It is not difficult to prove the Bollob\'as-Gy\'arf\'as Conjecture when $k = 1$ or $2$ (see \cite{BG}). We have extended this to the case $k = 3$.

\begin{thm}
For $n \ge 9$, $m(n,2,1,3) = n - 4$.
\end{thm}

The proof of this result involves a somewhat lengthy and delicate case analysis. We provide only a brief sketch, and refer the interested reader to~\cite{213} for a complete proof.

For $i,j \in \N$, define $K_{i,j}^G(i+j)$ to be a complete bipartite graph $K_{i,j} \subset G$. We simply write $K_{i,j}^G$ if $i$ and $j$ are known. Notice if $\overline{G}$ is not $k$-connected, then there exists a $K_{i,j}^G(|G|-k+1)$.

\begin{proof}
Let $f$ be a $2$-colouring of $E(K_n)$, and suppose that there is no monochromatic $3$-connected subgraph of $K_n$ on at least $n-4$ vertices. We shall show that there is a vertex of high degree in $R$ and in $B$. Since $R$ is not $3$-connected, there exists a $K_{i,j}^B(n-2)$. If $i \ge 3$, then this $K_{i,j}^B(n-2)$ is $3$-connected. If $i = 2$, then $j = n - 4$; let $I$ and $J$ be the partite sets of sizes $i$ and $j$, respectively. If $R[J]$ is $3$-connected, then it is the desired subgraph. Otherwise, $B[J]$ has a connected subgraph on $n - 6$ vertices, which, along with the vertices of $I$, form a $3$-connected subgraph on $n - 4$ vertices. Hence $i = 1$, so there is a vertex $x$ with $d_B(x) \ge n - 3$. Similarly, there is a vertex $y$ with $d_R(y) \ge n - 3$.

Assume without loss of generality that $xy \in E(R)$, and let $N \subset \Gamma_B(x)$ with $|N| = n - 3$. Writing $d_G^N(v)$ for $|\Gamma_G(v) \cap N|$, we can assume that $d^N_B(v) < n - 4$ for all $v \in N$, since otherwise there is a $K_{2,n-4}^B$ and we are done as before.

The remainder of the proof is an analysis of the following cases: either there is a $v \in N$ with $d_B^N(v) = n - 5$, or $d_R^N(v) \ge 2$ for all $v \in N$. In the latter case, we consider the three subcases corresponding to when the set $S = \{v \in N : d_R^N(v) = 2\}$ has cardinality $0$, $1$, or at least $2$.
\end{proof}

Bollob{\'a}s and Gy{\'a}rf{\'a}s noted that it is not even clear that in any $2$-colouring of $E(K_{4k-3})$, there is a monochromatic $k$-connected subgraph at all. A proof of this could probably be used to improve the bound $n \ge \min( (9+\sqrt{10})k, 13k-15)$ in Theorem~\ref{21k}.

\section{General $r$ and $s = 1$}\label{sectr1k}

In this section we consider the case $s=1$, but for general $r$. The question of Bollob\'as thus becomes, what is $m(n,r,1,k)$? We can derive an upper bound by considering finite affine planes.

\begin{lemma}\label{r1kupper}
Let $n,r,k \in \N$, with $n \ge r(k - 1)$ and $r - 1$ a prime
power. Then
$$m(n,r,1,k) \; \le \; \frac{n-k+1}{r-1} \: + \: r,$$ and if
$(r-1)^2$ divides $(n - r(k-1))$, then $m(n,r,1,k) \le \frac{n - k
+ 1}{r - 1}$.

Moreover $m(n,3,1,k) \le (n - k + 3)/2$ for every $n,k \in \N$,
and if $n \le 2r(k - 1)$, then $m(n,r,1,k) = 0$.
\end{lemma}

\begin{proof}
Let $n,r,k \in \N$, with $n \ge r(k - 1)$ and $r - 1$ a prime power. We shall describe a colouring $f$ of the edges of $K_n$ in which there is no monochromatic $k$-connected subgraph on more than $\frac{n-k+1}{r-1}$ vertices.

Since $r - 1$ is a prime power, there exists a finite affine plane $AF_{r-1}$ of order $r - 1$. Let $p_1, \ldots, p_{(r-1)^2}$ be the points and $P_1, \ldots, P_r$ be the parallel classes of $AF_{r-1}$. Let $C_1, \ldots, C_r$ be disjoint subsets of $V(K_n)$, each of size $k - 1$, and let $W = V(K_n) \setminus \bigcup_{i=1}^r C_i$. Now divide $W$ into $(r-1)^2$ classes $V_1, \ldots, V_{(r-1)^2}$ of about equal size (i.e., $|(|V_i| - |V_j|)| \le 1$ for every pair $i,j$).

The colouring $f$ is defined as follows. If $x \in V_i$ and $y \in
V_j$ are vertices of $K_n$ and $i \neq j$, then let $f(xy) = t$ if
and only if $p_i$ and $p_j$ lie on the same line in the class
$P_t$. If $i = j$ then $f(xy)$ may be chosen arbitrarily. If $x
\in C_i$, and $y \in W$, then let $f(xy) = i$. If $x \in C_i$ and
$y \in C_j$, then let $f(xy) = \min(i,j)$.

Let $\ell \in [r]$, and let $G$ be a monochromatic, $k$-connected
subgraph of $K_n$, with all edges coloured $\ell$ by $f$. Suppose
$G$ contains vertices from two different lines of $P_\ell$. Then
removing the vertices $V(G) \cap C_\ell$ from $G$ disconnects $G$,
and $|V(G) \cap C_\ell| \le k-1$, a contradiction. So $G$ contains
vertices from at most $r-1$ of the sets $V_i$. Similarly, $G$ may
contain no vertex of the set $C_i$ if $i \neq \ell$. Hence
\begin{equation}\label{r1kupp} |G| \; \le \; (r-1) \left\lceil
\frac{n-r(k-1)}{(r-1)^2} \right \rceil + k - 1 \; \le \;
\frac{n-k+1}{r-1} \: + \: r.\end{equation} Since $\ell$ and $G$
were arbitrary, this completes the proof of the first inequality,
and if $(r-1)^2$ divides $n - r(k-1)$, then we can remove the $r$
term from the right-hand side of (\ref{r1kupp}).

If $r = 3$, we split into two cases: $n - 3k + 3 \equiv 1 \pmod
4$, and $n - 3k + 3 \not\equiv 1 \pmod 4$. If $n - 3k + 3 = 4q +
1$, then exactly one of the sets $V_i$ has order $q + 1$, and so
(\ref{r1kupp}) becomes $|G| \le 2q + 1 + (k - 1) = (n - k + 2)/2$.
If $n - 3k + 3 \equiv 0, 2$ or $3 \pmod 4$, then $\left\lceil
\frac{n - 3k + 3}{4} \right \rceil \le \frac{n - 3k + 5}{4}$, so
$|G| \le \frac{n - 3k + 5}{2} + k - 1 = \frac{n - k + 3}{2}$.

To prove the final part of the lemma, let $n \le 2r(k-1)$, and
consider the following colouring $g$ of $E(K_n)$. First, partition
the vertices of $K_n$ into $2r$ sets $D_1, \ldots, D_{2r}$, each
of size at most $k - 1$. It is well-known (and easy to prove,
see~\cite{MGT} for example) that one can partition the edges of
$K_{2r}$ into $r$ edge-disjoint Hamilton paths of length $2r - 1$,
with each vertex an end-vertex of exactly one path; let these
paths be $Q_1, \ldots, Q_r$. If $x \in D_i$ and $y \in D_j$ with
$i \neq j$ and $ij \in Q_t$, then let $g(xy) = t$; if $i = j$, and
$i$ is an end-vertex of $Q_{t'}$, then let $g(xy) = t'$. It is
easy to check that the above colouring contains no $k$-connected
monochromatic subgraph, so if $n \le 2r(k-1)$ then $m(n,r,1,k) =
0$.
\end{proof}

Below is the colour $2$ subgraph of the colouring described in
Lemma~\ref{r1kupper} when $r = 3$ (Figure 2).

\[ \unit = 1.1cm
\medline \ellipse{0}{0.5}{1}{1} \ellipse{2.5}{0.5}{1}{1}
\ellipse{0}{3}{1}{1} \ellipse{2.5}{3}{1}{1}
\ellipse{-2.5}{0}{0.5}{0.5} \ellipse{-2.5}{1.75}{0.5}{0.5}
\ellipse{-2.5}{3.5}{0.5}{0.5} \thkline
\bez{-2.3}{1.4}{-0.5}{-2}{2}{-0.2}
\bez{-2.3}{2.1}{-0.5}{5.5}{2}{3.7} \dl{-2.2}{1.6}{-0.8}{0.8}
\dl{-2.2}{1.9}{-0.8}{2.7} \dl{0}{1.3}{0}{2.2}
\dl{2.5}{1.3}{2.5}{2.2} \dl{-2.5}{0.35}{-2.5}{1.35}
\ptlu{-2.5}{-0.2}{k - 1} \ptlu{-2.5}{1.55}{k - 1}
\ptlu{-2.5}{3.3}{k - 1} \ptll{-3.1}{3.5}{C_1}
\ptll{-3.1}{1.75}{C_2} \ptll{-3.1}{0}{C_3} \ptll{-1.05}{0.5}{V_3}
\ptll{-1.05}{3}{V_1} \ptlr{3.55}{0.5}{V_4} \ptlr{3.55}{3}{V_2}
\point{-0.7}{-1.7}{\small Figure 2}
\]

Lemma~\ref{r1kupper} gives us the upper bounds in Theorems~\ref{r1k} and \ref{31k}. We shall now show that $m(n,r,1,1) \ge \frac{n}{r-1}$ for every $n$ and $r$. F{\"u}redi~\cite{ZF} and Gy{\'a}rf{\'a}s~\cite{AG} discovered this while studying a hypergraph covering problem, namely, if one has $r$ partitions of $[n]$ such that every $x,y \in [n]$ lie in a common block of at least one of them, then how small can the largest block be? This is obviously equivalent to our problem, since the monochromatic components define $r$ partitions of $V(K_n)$, and if an edge is coloured $i$ then its endpoints lie in the same block of the $i^{th}$ partition.

We present a short, simple proof of this result, the ideas of which will be extended to give the lower bound in Theorem~\ref{r1k}.

\begin{lemma}\label{bipdensity}
Let $m,n \in \N$ and $c \in [0,1]$. If $G$ is a bipartite graph with part-sizes $m$ and $n$, and $e(G) \ge cmn$, then $G$ has a component of order at least $c(m+n)$.
\end{lemma}

\begin{proof}
If $c = 0$ the result is trivial, so assume $c > 0$. Let $M$ and $N$ be the partite sets of sizes $m$ and $n$, respectively, and let $xy \in E(G)$. The order of the component of $G$ containing $xy$ is at least $d(x) + d(y)$. Since $$\sum_{xy\in E(G)} (d(x) + d(y)) \; = \; \sum_{v \in V(G)} d(v)^2 \; = \; \sum_{v\in M} d(v)^2 \: + \: \sum_{v \in N} d(v)^2$$ $$\ge \; \left(\frac{e(G)}{m}\right)^2m \: + \: \left(\frac{e(G)}{n} \right)^2 n \; = \; \frac{e(G)^2(m+n)}{mn},$$ there must be an edge $xy$ with $d(x) + d(y) \ge \frac{e(G)(m+n)}{mn} \ge c(m+n)$. The order of the component of $G$ containing $xy$ is therefore at least $c(m+n)$.
\end{proof}

\begin{cor}\label{r11bip}
The order of the largest monochromatic component of an $r$-colouring of $E(K_{m,n})$ is at least $\frac{m+n}{r}$.
\end{cor}

This result is best possible, since if the partite sets are $M$ and $N$, and $|M|$ and $|N|$ are both divisible by $r$, then we may partition $M$ into parts $M_1, \ldots, M_r$ and $N$ into parts $N_1, \ldots, N_r$ of equal size, and colour all edges between $M_i$ and $N_j$ with colour $i-j \pmod r$. The largest monochromatic component in this colouring has order $(m+n)/r$.

\begin{thm}\label{r11}
Let $n,r \in \N$. Then $m(n,r,1,1) \ge \frac{n}{r-1}$.
\end{thm}

\begin{proof}
Let $n,r \in \N$, let $f$ be an $r$-colouring of $E(K_n)$, and let $C$ be a monochromatic component of $K_n$. If $C$ spans the whole of $V(K_n)$, then $M(f,n,r,1,k) = n$, and we are done. Otherwise, the edges of $K_n[C,V(K_n) \setminus C]$ are $(r-1)$-coloured by $f$, since $C$ is a (maximal) component. Thus, by Corollary~\ref{r11bip}, $K_n$ contains a monochromatic component of order at least $\frac{n}{r-1}$.
\end{proof}

We now return to the situation for general $k$. The strategy we shall use to prove the lower bound in Theorem~\ref{r1k} is analogous to that used in the proof of Theorem~\ref{r11}. First, in Lemma~\ref{r1kbip}, we shall derive an (asymptotically tight) upper bound on the number of edges in a bipartite graph with no large $k$-connected subgraph (as we did in Lemma~\ref{bipdensity}). From there we simply determine how large a $k$-connected subgraph this ensures.

We shall use the following simple observation in the proof of Lemma~\ref{r1kbip}.

\begin{lemma}{\label{arithmetic}}
If $a,b,c,d > 0$, then
$$\frac{ab}{a+b}+\frac{cd}{c+d} \; \le \; \frac{(a+c)(b+d)}{a+b+c+d}.$$
\end{lemma}

\begin{proof}
Expanding the inequality shows it is equivalent to $(ad-bc)^2\ge 0$.
\end{proof}

The next lemma is the key step in the proof of Theorem~\ref{r1k}. It is the analogue of Lemma~\ref{bipdensity} for general $k$.

\begin{lemma}{\label{r1kbip}}
Let $q,\ell,m,n \in \N$ with $m,n \ge \ell$ and $m + n \ge 2\ell + 1$. Let $G$ be a bipartite graph with parts $M$ and $N$ of size $m$ and $n$, respectively. If $G$ has no $(\ell+1)$-connected subgraph on at least $q$ vertices, then
\begin{eqnarray}\label{turan}
e(G) \; \le \; \frac{q(n - \ell)(m - \ell)}{m + n - 2\ell} \: + \:
(\ell^2 + \ell)(m + n - 2\ell).
\end{eqnarray}
\end{lemma}

\begin{proof}
We prove this by induction on $m + n$. To prove the base case, suppose that $m = \ell$. The inequality
$$\frac{q(n - \ell)(m - \ell)}{m + n - 2\ell} \: + \: (\ell^2 + \ell) (m + n - 2\ell) \: \ge \: mn,$$
reduces to $(\ell^2 + \ell)(n - \ell) \ge \ell n$, which holds if $n \ge \ell + 1$. Similarly this inequality is true if $n = \ell$ and $m \ge \ell + 1$. Since $e(G) \le e(K_{m,n}) = mn$, inequality (\ref{turan}) holds when $m + n = 2\ell + 1$.

So let $q,\ell,m,n \in \N$, $m,n \ge \ell + 1$, and assume that the statement of the lemma holds if $|M| + |N| \le m + n - 1$. Let $G$ be a bipartite graph, with parts $M$ and $N$ of size $m$ and $n$ respectively, and with no $(\ell+1)$-connected subgraph on at least $q$ vertices. Suppose first that $q \ge m + n + 1$. Then
$$\frac{q(n - \ell)(m - \ell)}{m + n - 2\ell} \; + \; (\ell^2 + \ell) (m + n - 2\ell) \; > \; (n - \ell)(m - \ell) \; + \; (\ell^2 + \ell)(m + n - 2\ell)$$ $$ \; = \; mn \; + \; \ell^2(m + n + 1 - 2\ell - 2) \; \ge \; mn \; \ge \; e(G),$$ and so inequality~(\ref{turan}) holds in this case.

Next suppose that $q \le m + n$. Since $G$ contains no $(\ell+1)$-connected subgraph on at least $q \le |G|$ vertices, $G$ itself cannot be $(\ell+1)$-connected, so there exists a cutset $C$ of size at most $\ell$. Let $x \in M$ and $y \in N$ be disconnected by $C$ (i.e. they are in different components of $G - C$). Since $m,n \ge \ell + 1$, we can choose a set $C' \supset C$, $x,y \notin C'$, with $|C' \cap M| = |C' \cap N|=\ell$. Since $x$ and $y$ were in different components of $G - C$, they must be in different components of its subgraph $G - C'$, so $G-C'$ is disconnected.

Let $G_1$ be a component of $G - C'$ and let $G_2 = G - (V(G_1) \cup C')$. For $i = 1,2$, let $H_i$ be the subgraph induced by $V(G_i) \cup C'$, and let $m_i = |V(H_i) \cap M|$ and $n_i = |V(H_i) \cap N|$. Note that since $|C' \cap N| = |C' \cap M| = \ell$, we have $m_i,n_i \ge \ell$, and $2\ell + 1 \le m_i + n_i \le m + n - 1$, since $V(G_1)$ and $V(G_2)$ are non-empty. Hence we can apply the induction hypothesis to the graphs $H_1$ and $H_2$, since if $H_i$ contains an $(\ell+1)$-connected subgraph on at least $q$ vertices then so does $G$.

Now $E(G) = E(H_1) \cup E(H_2)$, so $e(G) \le e(H_1) + e(H_2)$,
and by the induction hypothesis we have
\begin{eqnarray*}
e(H_1) \: + \: e(H_2) & \le & q \left( \frac{ (n_1 - \ell)(m_1 -
\ell )}{m_1 + n_1 - 2\ell } \: + \: \frac{(n_2 - \ell )(m_2 - \ell
)}{m_2 + n_2 - 2\ell} \right) \\[+1ex]
& & + \;\: (\ell^2 \: + \: \ell)(m_1 \: + \: m_2 \: + \: n_1 \: +
\: n_2 \: - \: 4\ell).
\end{eqnarray*}

Applying Lemma~\ref{arithmetic} with $a = n_1 - \ell$, $b = m_1 -
\ell$, $c = n_2 - \ell$ and $d = m_2 - \ell$, and using the
identities $m_1 + m_2 = m + \ell$ and $n_1 + n_2 = n + \ell$, we
have
\begin{eqnarray*}
\frac{(n_1 - \ell)(m_1 - \ell )}{m_1 + n_1 - 2\ell } \; + \;
\frac{(n_2 - \ell )(m_2 - \ell )}{m_2 + n_2 - 2\ell} \; \: \le
\\[+1ex]
\hspace{2.4cm} \le \;\: \frac{(n_1 + n_2 - 2\ell)(m_1 + m_2 -
2\ell)}{m_1 + m_2 + n_1 + n_2 - 4\ell} & = & \frac{(n - \ell)(m -
\ell)}{m + n - 2\ell},
\end{eqnarray*}
and hence
$$e(G) \; \le \; e(H_1) + e(H_2) \; \le \; q\frac{(n - \ell)(m -
\ell)}{m + n - 2\ell} \: + \: (\ell^2 + \ell)(m + n - 2\ell),$$ so
the induction step is complete. The lemma follows immediately.
\end{proof}

The lower bound in Theorem~\ref{r1k} now follows from
Lemma~\ref{r1kbip} and the following well-known theorem of
Mader~\cite{Mader}.

\begin{mader}
Let $\alpha \in \R$, and $G$ be a graph with average degree
$\alpha$. Then $G$ has an $\alpha/4$-connected subgraph.
\end{mader}

Note that since, in any $r$-colouring of $K_n$, some colour occurs at least $n(n-1)/2r$ times, Mader's Theorem implies the existence of a monochromatic $(n-1)/4r$-connected subgraph. This subgraph is $k$-connected if $n \ge 4kr + 1$, and has at least $(n-1)/4r + 1$ vertices. It is this weak bound that we shall need to prove the lower bound in Theorem~\ref{r1k}.

\begin{proof}[Proof of Theorem~\ref{r1k}]
Let $n,k,r \in \N$ with $k \ge 2$, $r \ge 3$ and $r - 1$ a prime power. The upper bound on $m(n,r,1,k)$ follows from Lemma~\ref{r1kupper}, so only the lower bound remains to be shown. If $n \le 11(k^2 - k)(r^2 - r)$ then the result holds vacuously, so assume $n > 11(k^2 - k)(r^2 - r)$. Let $f$ be an $r$-colouring of $E(K_n)$, and for $1 \le i \le r$ let $G^{(i)}$ denote the graph on $V(K_n)$ with edge set $f^{-1}(i)$. We shall find, for some $i \in [r]$, a $k$-connected subgraph of $G^{(i)}$ on at least $\frac{n}{r-1} - 11(k^2-k)r$ vertices.

Let $H$ be a monochromatic $k$-connected subgraph of $K_n$ of maximum order, and suppose without loss that $H$ has colour $1$. Let $C = V(H)$, $|C| = c$, $D = V(K_n) \setminus C$ and $|D| = d$. By Mader's Theorem, $c \ge (n-1)/4r + 1 > n/4r$, and we may assume that $c < n/(r-1)$, since otherwise $H$ is the desired monochromatic subgraph. Thus $c,d > k$, since $r \ge 3$ and $n > 4kr$. We shall apply Lemma~\ref{r1kbip} to the bipartite graph $G^{(i)}[C,D]$, where $i \in [2,r]$ is chosen to maximize the number of edges in this graph.

Since $H$ is maximal, no vertex of $D$ sends more than $k-1$ edges of colour $1$ into $C = V(H)$, so by the pigeonhole principle, for some $i \in [2,r]$ there are at least $d(c-k+1)/(r-1)$ edges between $C$ and $D$ of colour $i$. Fix this $i$, let $\ell = k - 1$ and let $G = G^{(i)}[C,D]$. By Lemma~\ref{r1kbip}, if $q \in
\N$ satisfies
$$q \frac{(d - \ell)(c - \ell)}{(c + d - 2\ell)} \: + \: (\ell^2 + \ell)
(c + d - 2\ell) \; < \; \frac{d(c - \ell)}{r - 1} \; \le \;
e(G),$$ or, equivalently,
\begin{eqnarray}\label{r1kpf}
q \; < \; \frac{d(c + d - 2\ell)}{(d - \ell)(r - 1)} \: - \:
\frac{(\ell^2 + \ell)(c + d - 2\ell)^2}{(c - \ell)(d - \ell)},
\end{eqnarray} then $G$ contains a $k$--connected subgraph on at
least $q$ vertices.

The theorem will now follow if we can show that the right-hand side of (\ref{r1kpf}) is greater than $\frac{n}{r-1} - 11(k^2 - k)r$, by setting $q$ equal to this value. Since $d > d - \ell$, $c + d = n$ and $\ell^2 + \ell = (k - 1)^2 + (k - 1) = k^2 - k$, we have $$\frac{d(c + d - 2\ell)}{(d - \ell)(r - 1)} \; - \; \frac{(\ell^2 + \ell)(c + d - 2\ell)^2}{(c - \ell)(d - \ell)} \; > \; \frac{n - 2\ell}{r - 1} \; - \; \frac{(k^2 - k)(n - 2\ell)^2}{(c - \ell)(d - \ell)}.$$ It therefore only remains to bound $(c - \ell)(d - \ell)$ from below. Since $c + d = n$, $(c - \ell)(d - \ell)$ is increasing with $c$ for $c < n/2$, so since $c < n/(r-1)$ and $r \ge 3$, the minimum is achieved by taking $c$ to be as small as possible. Hence, by setting $c = n/4r$, we get
\begin{eqnarray*}
(c - \ell)(d - \ell) & > & \left( \frac{n}{4r} - \ell \right) \left( \left( \frac{4r - 1}{4r} \right) n - \ell \right)\\ & > & \left( \frac{n}{8r} \right) \left( \frac{(8r - 3)n}{8r} \right) \; > \; \frac{n^2}{10r}, \end{eqnarray*}
since $n > 11(k^2 - k)(r^2 - r) > 8\ell r$ and $r \ge 3$. We have therefore shown that if $$q \; = \; \left\lfloor \frac{n - 2\ell}{r - 1} \: - \: \frac{(k^2 - k)n^2} {(n^2/10r)} \right\rfloor \; > \; \frac{n}{r - 1} \: - \: 11(k^2 - k)r,$$ then by Lemma~\ref{r1kbip} there exists a monochromatic $k$-connected subgraph on at least $q$ vertices. This completes the proof.
\end{proof}

Having proved Theorem~\ref{r1k}, we can now use it (in place of Mader's Theorem) to give the following slight improvement for sufficiently large values of $n$.

\begin{thm}\label{r1kextra}
Let $n,k,r \in \N$ and $\eps > 0$ satisfy $r \ge 3$ and $n \ge \frac{11(2 + \eps)}{\eps}k^2r^2$. Then
$$m(n,r,1,k) \; \ge \; \frac{n}{r-1} \: - \: \left(1 + \frac{1}{r(r-2)} + \eps \right)k^2r.$$
In particular, if $n \ge 44k^2r^2$, then $m(n,r,1,k) \ge \frac{n}{r-1} - 2k^2r$.
\end{thm}

\begin{proof}
The proof follows exactly as the proof of Theorem~\ref{r1k}, but we can now give the following improved bound on $(c - \ell)(d - \ell)$, since we know $c > \frac{n}{r-1} - 11k^2r$.
\begin{eqnarray*}
(c - \ell)(d - \ell) & > & \left( \frac{n}{r-1} \: - \: 11k^2r \: - \: \ell \right) \left( \frac{(r-2)n}{r-1} \: + \: 11k^2r \: - \: \ell \right)\\[+1ex]
& > & \left( \frac{r - 2}{(r - 1)^2} \right)n^2 \: - \: \left( 11k^2 \left( \frac{r^2 - 3r}{r-1} \right) + \ell \right)n - 121k^4r^2\\[+1ex] & > & \left( \frac{r - 2}{(r - 1)^2} \right)n^2 \: - \: 11k^2(r-1)n
\end{eqnarray*}
since $n \ge 13k^2r^2$ if $\eps \le 11$. Let $\delta = \eps (r-2)(1 + \frac{1}{r(r-2)} + \eps)^{-1}$. Then $\delta > \eps(r - 2)/(2 + \eps)$, so $$\delta n \; > \; 11k^2r^2(r-2) \; > \; 11k^2(r - 1)^3,$$ since $n \ge (11(2 + \eps)k^2r^2)/\eps$ and $r \ge 3$. Thus
\begin{eqnarray*}
(c - \ell)(d - \ell) & > & \left( \frac{r - 2 - \delta}{(r - 1)^2} \right)n^2,
\end{eqnarray*} so if
\begin{eqnarray*}
q \; = \; \frac{n-2\ell}{r-1} \: - \: \frac{(k^2 - k)(r - 1)^2}{r - 2 - \delta} \; > \; \frac{n}{r-1} \: - \:  \frac{k^2(r-1)^2}{r - 2 - \delta},
\end{eqnarray*}
then there exists a monochromatic $k$-connected subgraph on at least $q$ vertices. Now simply observe that we chose $\delta$ so that $$\frac{(r-1)^2}{r - 2 - \delta} = \left( 1 + \frac{1}{r(r-2)} + \eps \right)r,$$ and the theorem follows. The final implication is attained by setting $\eps = 2/3$ and recalling that $r \ge 3$.
\end{proof}

It would be interesting to know where in the ranges given by Theorems~\ref{r1k} and \ref{r1kextra} the truth lies. We strongly suspect that the upper bound from Lemma~\ref{r1kupper} gives the
correct answer.

\begin{conjecture}\label{conjr1k}
Let $n,k,r \in \N$ with $r \ge 3$, $n \ge 2r(k-1) + 1$, $r-1$ a prime power and $n - r(k-1)$ divisible by $(r-1)^2$. Then $$m(n,r,1,k) \; = \; \frac{n-k+1}{r-1}.$$
\end{conjecture}

\begin{rmk}\label{r1krmk}
By Lemma~\ref{r1kupper}, $m(n,r,1,k) = 0$ if $n \le 2r(k-1)$. Hence the lower bound on $n$ in the conjecture cannot be weakened any further.
\end{rmk}

We also have the following conjecture for the bipartite version of the question. It says that the order of the largest $k$-connected subgraph equals the upper bound given in Corollary~\ref{r11bip} (and so does not depend on $k$), as long as the partite sets are large.

\begin{conjecture}\label{conjr1kbip}
Let $m,n,k,r \in \N$, with $r \ge 3$ and $m,n \ge rk$. Any $r$-colouring of the edges of $K_{m,n}$ contains a monochromatic $k$-connected subgraph on at least $\frac{m+n}{r}$ vertices.
\end{conjecture}

Although we have been unable to prove Conjectures~\ref{conjr1k} and \ref{conjr1kbip}, Theorem~\ref{31k} shows that Conjecture~\ref{conjr1k} holds in the case $r = 3$. We shall next prove this result. We begin with an easy lemma.

\begin{lemma}\label{31kbip}
Let $k,p,q \in \N$ satisfy $3p \ge q \ge p \ge 24k$, and let $P$ and $Q$ be sets with $|P| = p$ and $|Q| = q$. Let $K_{p,q}$ be the complete bipartite graph with parts $P$ and $Q$. Suppose the edges of $K_{p,q}$ are $3$-coloured in such a way that each vertex in $P$ sends at most $k$ edges of colour $3$ into $Q$, and each vertex in $Q$ sends at most $k$ edges of colour $2$ into $P$.

Then the subgraph induced by edges of colour $1$ contains a $k$-connected subgraph $G$ with $|P \setminus V(G)| \le 16k$, and $|Q \setminus V(G)| \le 8k$. In particular, $|V(G)| \ge p + q - 24k$.
\end{lemma}

\begin{proof}
Let $k,p,q \in \N$ satisfy $3p \ge q \ge p \ge 24k$, and let $f$ be a 3-colouring of $E(K_{p,q})$ satisfying the conditions of the lemma. Let $$S_P = \{v \in P : v\textup{ sends at most }3q/4\textup{ edges of colour }1\textup{ into }Q\}\textup{, and}$$ $$S_Q = \{v \in P : v \textup{ sends at most }3p/4 \textup{ edges of colour }1\textup{ into }P\} \hspace{0.95cm}$$ be sets of `bad' vertices. We shall remove the bad sets and apply Lemma~\ref{intersect}.

We need to bound $|S_P|$ and $|S_Q|$ from above. Since each vertex of $Q$ has at most $k$ incident edges of colour $2$, we have $|f^{-1}(2)| \le kq$, and similarly $|f^{-1}(3)| \le kp$. Also, since each vertex of $S_P$ has at least $q/4$ incident edges of colour $2$ or $3$, we have $|f^{-1}(2)| + |f^{-1}(3)| \ge |S_P|(q/4)$, and similarly $|f^{-1}(2)| + |f^{-1}(3)| \ge |S_Q|(p/4)$. Thus
$$|S_P| \; \le \; \frac{4}{q}\left( |f^{-1}(2)| + |f^{-1}(3)| \right) \; \le \; \frac{4k(p+q)}{q} \; \le \; 8k, \: \textup{ and}$$
$$|S_Q| \; \le \; \frac{4}{p}\left( |f^{-1}(2)| + |f^{-1}(3)| \right) \; \le \; \frac{4k(p+q)}{p} \; \le \; 16k. \hspace{0.8cm}$$

Now, let $P' = P \setminus S_P$ and $Q' = Q \setminus S_Q$, and let $G$ be the bipartite graph with vertex set $P' \cup Q'$, and edge set $f^{-1}(1)$. If $x \in P'$, then $x$ sends at least $3q/4$ edges of colour $1$ into $Q$, so $$d_G(x) \; \ge \; 3q/4 - |S_Q| \; \ge \; 18k - 16k \; > \; k,$$ and similarly if $y,z \in Q'$, then
\begin{eqnarray*}
|\Gamma_G(y) \cap \Gamma_G(z)| & \ge & 3p/4 \: + \: 3p/4 \: - \: p
\: - \: |S_P|\\[+1ex]
& = & p/2 - |S_P| \; \ge \; 12k - 8k \; > \; k,
\end{eqnarray*} so the conditions of Lemma~\ref{intersect} are satisfied. Thus by Lemma~\ref{intersect}, $G$ is $k$-connected. Since also $|P \setminus V(G)| = |S_P| \le 8k$ and $|Q \setminus V(G)| = |S_Q| \le 16k$, $G$ is the desired subgraph.
\end{proof}

Given graphs $G$ and $H$, define $G \cup H$ to be the graph with vertex set $V(G) \cup V(H)$ and edge set $E(G) \cup E(H)$. We shall also use the following trivial observation.

\begin{obs}\label{GcupH}
Let $k \in \N$. If $G$ and $H$ are $k$-connected graphs, and $|V(G) \cap V(H)| \ge k$, then the graph $G \cup H$ is also $k$-connected.
\end{obs}

We are now ready to prove Theorem~\ref{31k}.

\begin{proof}[Proof of Theorem~\ref{31k}]
Let $n,k \in \N$ with $n \ge 480k$. The upper bound on
$m(n,3,1,k)$ follows from Lemma~\ref{r1kupper}, so only the lower
bound remains to be shown.

Let $f$ be a $3$-colouring of the edges of $K_n$, and let $V = V(K_n)$. For each $i \in \{1,2,3\}$, let $G_i$ be the subgraph of $K_n$ with vertex set $V$ and edge set of $f^{-1}(i)$, and assume that $G_i$ has no $k$-connected subgraph on more than $(n-k)/2$ vertices. We begin by covering $V$ with monochromatic
$k$-connected
subgraphs.\\[+2ex]
\underline{Claim 1}: There exist (not necessarily disjoint) subsets $A_1$, $A_2$, and $A_3$ of $V$ such that $G_i[A_i]$ is $k$-connected, and $A_1 \cup A_2 \cup A_3 = V$.

\begin{proof}
Assume, without loss of generality, that $e(G_1) \ge e(G_2) \ge e(G_3)$. By Mader's Theorem (and since $n > 12k$), there exists a maximal set $A_1 \subset V$, with $|A_1| \ge n/12$, such that $G_1[A_1]$ is $k$-connected. If $|A_1| > (n-k)/2$, then $G_1[A_1]$ is a monochromatic $k$-connected subgraph on more than $(n-k)/2$ vertices, contradicting our assumption, so (writing $A^c_1$ for $V \setminus A_1$) we have $|A^c_1| > n/2 > |A_1|$.

For $i = 2,3$, let $H_i = G_i[A_1,A^c_1]$ be the bipartite graph induced by the edges of colour $i$ and the sets $A_1$ and $A^c_1$. Since $A_1$ is maximal, each vertex of $A^c_1$ sends at most $k-1$ edges of colour $1$ into $A_1$ (by Observation~\ref{addvtx}), and so has degree at least $|A_1| - k + 1$ in $H_2 \cup H_3$. Hence
$$e(H_2) + e(H_3) \: \ge \: \left| A^c_1 \right| \left( |A_1| - k + 1 \right) \: \ge \: \frac{11n}{12} \left( \frac{n}{12} - k + 1 \right) \: > \: \frac{n^2}{15},$$ the second inequality holding because the function $-x^2 + (n - k + 1)x$ is increasing for $x < (n - k + 1)/2$, and the third holding because $n > 95k$.

Since $e(H_2) \ge e(H_3)$, we obtain $e(H_2) > n^2/30$, so the average degree in $H_2$ is at least $n/15$. Applying Mader's Theorem again, we deduce that $H_2$ contains an $(n/60)$-connected subgraph $H'_2$. Since $H_2$ is bipartite, $H'_2$ must contain at least $n/60$ vertices of each class of $H_2$; in particular, it must contain at least $8k$ vertices of $A_1$ (since $n \ge 480k$).

Let $A_2$ be a maximal set containing $V(H'_2)$ such that $G_2[A_2]$ is $k$-connected. We have now found sets $A_1$ and $A_2$, with $G_i[A_i]$ $k$-connected for $i = 1,2$. We complete the proof by using Lemma~\ref{intersect} to find a $k$-connected graph in $G_3$ containing $(A_1 \cup A_2)^c$.

Let $X = A_1 \cap A_2$ and $Y = (A_1 \cup A_2)^c$. Notice that $|A_2| \le (n-k)/2$, since otherwise we would have a monochromatic $k$-connected subgraph on more than $(n-k)/2$ vertices, contradicting our assumption. Since $V(H'_2) \subset A_2$ and, as observed above, $H'_2$ contains at least $8k$ vertices of $A_1$, we have $|X| \ge 8k$. Hence also
\begin{eqnarray*}
|Y| & = & n - |A_1| - |A_2| + |A_1 \cap A_2|\\[+1ex]
& \ge & n - (n - k) + |X| \; = \; |X| + k \; \ge \; 9k,
\end{eqnarray*}
since $|A_1|, |A_2| \le (n-k)/2$.

We want to apply Lemma~\ref{intersect} to the bipartite graph $G_3[X,Y]$, but first we must remove the vertices of degree at most $k-1$ from $X$, as in the proof of Theorem~\ref{21k}. As in that proof, let $$U = \{ v \in X : |\Gamma_{G_3}(v) \cap Y| \le k-1\}.$$

Since $G[A_1]$ and $G[A_2]$ are maximal monochromatic $k$-connected subgraphs, each vertex $v \in Y$ can send only at most $k-1$ edges of colour $1$ into $A_1$, and $k-1$ edges of colour $2$ into $A_2$. Therefore $v$ must send at least $|X| - 2k + 2$ edges of colour $3$ into $X = A_1 \cap A_2$. This is true for every $v \in Y$, so $G_3[X,Y]$ has at most $|Y|(2k-2)$ non-edges. But each vertex of $U$ sends at least $|Y|-k+1$ non-edges into $Y$. Thus we have
$$|U|(|Y|-k+1) \; \le \; |Y|(2k-2),$$
and hence
$$|U| \; \le \; \frac{2|Y|(k-1)}{|Y|-k+1} \; \le \; \frac{18k(k-1)}{8k+1} \; < \; 3k,$$ since the function
$\frac{2x(k-1)}{x-k+1}$ is decreasing for $x > k-1$, and $|Y| \ge
9k$.

Let $X' = X \setminus U$, and consider the bipartite graph $G_3[X',Y]$. By the definition of $U$, each vertex in $X'$ has degree at least $k$ in this graph. Also, as noted above, each vertex of $Y$ sends at most $2k-2$ edges of colour $1$ or $2$ into $X'$, so any two vertices in $Y$ have at least $$|X'| - 4k + 4 \; = \; |X| - |U| - 4k + 4 \; > \; 8k - 3k - 4k \; = \; k$$ common neighbors in $X'$.

So $G_3[X',Y]$ satisfies the conditions of Lemma~\ref{intersect}, and therefore by that lemma $G_3[X',Y]$ is $k$-connected. Let $A_3$ be a maximal set such that $G_3[A_3]$ is $k$-connected, and $X' \cup Y \subset A_3$. Since $V \setminus (A_1 \cup A_2) = Y \subset A_3$, this completes the proof of Claim 1.
\end{proof}

For the remainder of the proof, $\{i,j,\ell\}$ will always be the set $\{1,2,3\}$, though the order will vary. Let $A_1, A_2, A_3$ be the (maximal) sets given by Claim 1, and for each $i$ (i.e., for each triple $i,j,\ell$ with $\{i,j,\ell\} = \{1,2,3\}$), let $a_i = |A_i \setminus (A_j \cup A_\ell)|$, $b_i = |(A_j \cap A_\ell) \setminus A_i|$, and $c = |A_1 \cap A_2 \cap A_3|$ (see Figure 3).

\[ \unit = 1.1cm
\medline \ellipse{0}{0.7}{1.4}{1.4} \ellipse{-0.8}{2.1}{1.4}{1.4}
\ellipse{0.8}{2.1}{1.4}{1.4} \ptll{-1.7}{3.4}{A_1}
\ptlr{1.7}{3.4}{A_2} \ptlr{1}{-0.5}{A_3} \ptlu{-1.2}{2.3}{a_1}
\ptlu{1.2}{2.3}{a_2} \ptlu{0}{0}{a_3} \ptlu{0.8}{1}{b_1}
\ptlu{-0.8}{1}{b_2} \ptlu{0}{2.3}{b_3} \ptlu{0}{1.4}{c}
\point{-0.7}{-1.7}{\small Figure 3}
\]\\[-1ex]

By Claim 1,
\begin{equation}\label{31ka} \sum_i a_i + \sum_i b_i + c = n. \end{equation}
Our initial assumption says that $|A_i|= a_i + b_j + b_\ell + c \le (n-k)/2$ for each triple $i,j,\ell$. Summing over $i = 1,2,3$ and subtracting (\ref{31ka}) gives
\begin{equation}\label{31kb} \sum_i b_i + 2c \; \le \; \frac{n-3k}{2},\end{equation} whilst summing pairwise and subtracting (\ref{31ka}) gives
\begin{equation}\label{31kc} a_i \; \ge
\; b_i + c + k \end{equation} for each $i \in \{1,2,3\}$.

Now, observe that since $G_j[A_j]$ is a maximal monochromatic $k$-connected subgraph, each vertex of $A_i \setminus A_j$ sends at most $k-1$ edges of colour $j$ into $A_j \setminus A_i$, for each pair $i,j$. We wish to apply Lemma~\ref{31kbip} to the pair of sets $A_i \setminus A_j$ and $A_j \setminus A_i$; the next claim (which we shall also prove using Lemma~\ref{31kbip}) allows us to do so.\\[+2ex]
\underline{Claim 2}: $a_i \ge n/6$ for each $i \in \{1,2,3\}$.

\begin{proof}
Let $i \in \{1,2,3\}$ and suppose $a_i < n/6$. Note that by (\ref{31kc}) we also have $b_i + c \le a_i - k < n/6$. Assume, without loss of generality, that $|A_j \setminus A_\ell| \le |A_\ell \setminus A_j|$. We shall apply Lemma~\ref{31kbip} with $P = A_j \setminus A_\ell$ and $Q = A_\ell \setminus A_j$.

Let $p = |A_j \setminus A_\ell|$ and $q = |A_\ell \setminus A_j|$. By assumption, $q \ge p$. Now observe that $p \ge 24k$, since otherwise
\begin{eqnarray*}
|A_j| & = & n \: - \: |A_\ell \setminus A_j| \: - \: |A_i \setminus (A_j \cup A_\ell)|\\[+1ex]
& = & n - p - a_i \; > \; \frac{5n}{6} - 24k \; > \; \frac{n}{2},
\end{eqnarray*}
since $a_i < n/6$ and $n > 72k$, which contradicts our assumption that $|A_j| \le (n-k)/2$. Also note that $q \le 3p$, since
\begin{eqnarray*}
p \: + \: q & = & |A_j \bigtriangleup A_\ell| \; = \; n - |A_i\setminus (A_j \cup A_\ell)| \: - \: |A_j \cap A_\ell|\\ & = & n \: - \: (a_i + b_i + c),
\end{eqnarray*} so if $q > 3p$, then
\begin{eqnarray*} |A_j| & = & |A_j \setminus A_\ell| \: + \: |(A_j \cap A_\ell) \setminus A_i| \: + \: |A_j \cap A_\ell \cap A_i|\\ & = & q \: + \: b_i \: + \: c \;\: > \;\: \frac{3(n - a_i - b_i - c)}{4} \: + \: b_i \: + \: c\\
& > & \frac{3(n - a_i)}{4} \; > \; \frac{5n}{8},
\end{eqnarray*}
which again contradicts our assumption that $|A_j| \le (n-k)/2$.

Hence $k$, $p$ and $q$ satisfy $3p \ge q \ge p \ge 24k$. Also, as observed above, each vertex of $P = A_j \setminus A_\ell$ sends at most $k-1$ edges of colour $\ell$ into $Q =  A_\ell \setminus A_j$, and similarly each vertex of $Q$ sends at most $k-1$ edges of colour $j$ into $P$, by maximality of $A_\ell$ and $A_j$. So by Lemma~\ref{31kbip}, there must exist a monochromatic $k$-connected subgraph in $G_i[P,Q]$ on at least $p + q - 24k$ vertices. Since
$$p + q - 24k \; = \; n - (a_i + b_i + c) - 24k \; > \; \frac{2n}{3} - 24k \; > \; \frac{n}{2}$$(because $b_i + c < a_i < n/6$ and $n > 144k$), this contradicts our assumption that there is no monochromatic $k$-connected subgraph in $G_i$ on more than $(n-k)/2$ vertices. This final contradiction completes the proof
of the claim.
\end{proof}

We shall now apply Lemma~\ref{31kbip} to the sets $A_i \setminus A_j$ and $A_j \setminus A_i$, for each pair $i$ and $j$. Let $i,j \in \{1,2,3\}$ with $i \neq j$, and assume, without loss of generality, that $|A_i \setminus A_j| \le |A_j \setminus A_i|$. We shall apply Lemma~\ref{31kbip} with $P = A_i \setminus A_j$ and $Q = A_j \setminus A_i$.

Let $p = |A_i \setminus A_j|$ and $q = |A_j \setminus A_i|$. By assumption, $q \ge p$. Now, $p \ge |A_i \setminus (A_j \cup A_\ell)| = a_i$, so by Claim 2, $p \ge n/6$. Since $n \ge 144k$, it follows that $p \ge 24k$, and since $q \le |A_j| < n/2$, it also follows that $3p \ge q$. As observed earlier, each vertex of $P$ sends at most $k-1$ edges of colour $j$ into $Q$, and each vertex of $Q$ sends at most $k-1$ edges of colour $i$ into $P$, since $A_i$ and $A_j$ are maximal.

Therefore, applying Lemma~\ref{31kbip} to the sets $P$ and $Q$, we obtain a $k$-connected subgraph of $G_\ell$ (where, as usual, $\ell = \{1,2,3\} \setminus \{i,j\}$) omitting at most $16k$ vertices of $P$ and at most $8k$ vertices of $Q$. Let this subgraph be $L_\ell$.

We obtain in this way three $k$-connected subgraphs, $L_1$, $L_2$ and $L_3$. For each $\ell \in \{1,2,3\}$, let $M_\ell$ be the vertex set of a maximal $k$-connected subgraph of $G_\ell$ containing $L_\ell$. Now, for each pair $i \neq j$, let $X_{ij} = A_i \setminus (A_j \cup M_\ell)$ be the set of vertices in $A_i \setminus A_j$ avoided by $M_\ell$, and let $x_{ij} = |X_{ij}| \le 16k$. Also, for each $\ell \in \{1,2,3\}$, let $Z_\ell = M_\ell \setminus (A_i \cup A_j)$, and let $z_\ell = |Z_\ell|$. We have, therefore, for each triple $i,j,\ell$, that \begin{equation}\label{31kcc} a_i + a_j + b_i + b_j - x_{ij} - x_{ji} + z_\ell \: \le \: |M_\ell| \: \le \: \frac{n-k}{2}, \end{equation} by assumption, since $M_\ell$ is the vertex set of a monochromatic $k$-connected subgraph.

Although we have so far been approximating wildly, we must now be precise. Summing the inequalities (\ref{31kcc}) over $\ell = 1,2,3$, we get
$$2\sum_i a_i \: + \: 2\sum_i b_i \: - \: \sum_{i,j} x_{ij} \: + \: \sum_i z_i \; \le \; \frac{3(n-k)}{2},$$ which is equivalent to
$$\frac{n+3k}{2} \; \le \; 2c \: + \: \sum_{i,j} x_{ij} \: - \: \sum_i z_i,$$ since $\sum a_i + \sum b_i + c = n$. Combining this with $\sum b_i + 2c \le (n-3k)/2$, we obtain
\begin{equation}\label{31kd}
0 \; \le \; \sum_i b_i \; \le \; \frac{n-3k}{2} - 2c \; \le \;
\sum_{i,j} x_{ij} - \sum_i z_i - 3k,
\end{equation} so $\sum x_{ij} - \sum z_i \ge 3k$, and by the
pigeonhole principle, there exists an $i \in \{1,2,3\}$ such that
\begin{equation}\label{31ke} x_{ij} + x_{i\ell} - z_i \ge
k.\end{equation}\\[-1ex]
We fix this $i$ for the remainder of the proof. We shall show that
inequality (\ref{31ke}) implies that $|X_{ij} \cup X_{i\ell}| \ge
k$, and deduce that $G_i[A_i \cup M_i]$ is $k$-connected.

Indeed, let $j \in \{1,2,3\} \setminus \{i\}$, and consider a vertex $v \in X_{ij}$. We shall show that $v \in M_i$ and therefore that $X_{ij} \cap X_{i\ell} \subset Z_i$. Let $\ell \in \{1,2,3\} \setminus \{i,j\}$, and recall that $X_{ij} \subset A_i \setminus (A_j \cup M_\ell)$, so $v \notin A_j$ and $v \notin M_\ell$. Since $v \notin A_j$ and $A_j$ is maximal, $v$ sends at most $k-1$ edges of colour $j$ into $A_j$, and since $v \notin M_\ell$ and $M_\ell$ is maximal, $v$ sends at most $k-1$ edges of colour $\ell$ into $M_\ell$.

How many edges of colour $j$ or $\ell$ can $v$ send into $A_j \setminus (A_i \cup A_\ell)$? Since $M_\ell$ contains $V(L_\ell)$, we know that $M_\ell$ avoids at most $16k$ vertices of $A_j \setminus A_i$, so by the observations above, $v$ sends at most $17k - 1$ edges of colour $\ell$ into $A_j \setminus A_i$, and so at most $18k - 2$ edges of colour $j$ or $\ell$ into $A_j \setminus A_i$. Thus $v$ sends at most $18k - 2$ edges of colour $j$ or $\ell$ into $A_j \setminus (A_i \cup A_\ell) \subset A_j \setminus A_i$.

Now, $M_i$ avoids at most $16k$ vertices of the set $A_j \setminus A_\ell$, and so at most $16k$ vertices of $A_j \setminus (A_i \cup A_\ell)$. Therefore, the number of edges of colour $i$ going from $v$ into $M_i$ is at least
\begin{eqnarray*}
|A_j \setminus (A_i \cup A_\ell)| \: - \: 16k \: - \: (18k - 2) & = & a_i \: - \: 34k \: + \: 2 \\
\; > \hspace{0.3cm} \frac{n}{6} \: - \: 34k & > & k,
\end{eqnarray*} since $a_i \ge n/6$ by Claim 2, and $n > 210k$. But $M_i$ was chosen to be a maximal monochromatic $k$-connected subgraph, so if $v$ sends at least $k$ edges of colour $i$ into $M_i$, it follows that $v \in M_i$.

Now, suppose that in fact $v \in X_{ij} \cap X_{i\ell}$. Since $X_{ij} \subset A_i \setminus A_j$ and $X_{i\ell} \subset A_i \setminus A_\ell$, it follows that $v \in A_i \setminus (A_j \cup A_\ell)$. Recall that $Z_i = M_i \setminus (A_j \cup A_\ell)$ and it is clear that, in this case, $v \in M_i$ implies $v \in Z_i$. Hence $X_{ij} \cap X_{i\ell} \subset Z_i$, and so $z_i \ge |X_{ij} \cap X_{i\ell}|$.

It now follows immediately that $$|X_{ij} \cup X_{i\ell}| \; = \; |X_{ij}| + |X_{i\ell}| - |X_{ij} \cap X_{i\ell}| \; \ge \; x_{ij} + x_{i\ell} - z_i \; \ge \; k,$$ by inequality (\ref{31ke}). The following claim now gives us the final contradiction.\\[+2ex] \underline{Claim 3}: $G_i[A_i \cup M_i]$ is $k$-connected, and has order at least $3n/4$.

\begin{proof}
Let $i$ be as in inequality (\ref{31ke}), and $\{i,j,\ell\} = \{1,2,3\}$. We have shown that $|X_{ij} \cup X_{i\ell}| \ge k$, and that $X_{ij} \cup X_{i\ell} \subset A_i \cap M_i$. By the definitions of $A_i$ and $M_i$, the graphs $G_i[A_i]$ and $G_i[M_i]$ are $k$-connected. Therefore, by Observation~\ref{GcupH}, $G_i[A_i \cup M_i]$ is $k$-connected.

Now, since $M_i$ contains $V(L_i)$, we know that $M_i$ avoids at most $24k$ vertices of $A_j \bigtriangleup A_\ell$, so $A_i \cup M_i$ avoids at most $24k + |(A_j \cap A_\ell) \setminus A_i| = 24k + b_i$ vertices of $V$. By inequality (\ref{31kd}), we have (very weakly), that $b_i \le \sum_m b_m \le \sum_{u,v} x_{uv} \le 96k$, since $x_{uv} \le 16k$ for each $u \neq v$, $u,v \in \{1,2,3\}$. Thus $$|A_i \cup M_i| \; \ge \; n - (24k + b_i) \; \ge \; n - 120k \; \ge \; \frac{3n}{4},$$ since $n \ge 480k$. This completes the proof of the claim.
\end{proof}

So $G_i[A_i \cup M_i]$ is a monochromatic $k$-connected subgraph on more than $(n-k)/2$ vertices, contradicting our assumption that no such subgraph exists. This contradiction proves the theorem.
\end{proof}

It is now easy to obtain the exact value of $m(n,3,1,k)$ whenever $n \ge 480k$.

\begin{cor}\label{31kexact}
Let $n,k \in \N$, with $n \ge 480k$. Then\\
\begin{equation*}
m(n,3,1,k) =
\begin{cases}
(n - k + 1)/2 & \text{\emph{if } $n + k \equiv 1$ $(\emph{mod}$
$4)$,}\\
(n - k + 2)/2 & \text{\emph{if } $n + k \equiv 0$ \emph{or} $2$
$(\emph{mod}$ $4)$,}\\
(n - k + 3)/2 & \text{\emph{if } $n + k \equiv 3$ $(\emph{mod}$
$4)$.}
\end{cases}
\end{equation*}
\end{cor}

\begin{proof}
Let $n,k \in \N$, with $n \ge 480k$. If $n + k \equiv 1 \pmod 4$, then $n - 3k + 3 \equiv 0 \pmod 4$, so by Lemma~\ref{r1kupper} and Theorem~\ref{31k} we have $m(n,3,1,k) = (n - k + 1)/2$. If $n + k \equiv 0$ or $2 \pmod 4$, then $(n - k + 2)/2$ is the only integer in the range given by Theorem~\ref{31k}, so clearly $m(n,3,1,k) = (n - k + 2)/2$. If $n + k \equiv 3 \pmod 4$, then we have $m(n,3,1,k) \le (n - k + 3)/2$ by Lemma~\ref{r1kupper}.

It remains to prove the lower bound in the case $n + k \equiv 3 \pmod 4$. To do this, we follow the proof of Theorem~\ref{31k}, making a couple of small alterations.

To be precise, let $n,k \in \N$, with $n \ge 480k$, and $n + k \equiv 3 \pmod 4$. Let $f$ be a 3-colouring of $E(K_n)$, and assume that $G_i$ contains no monochromatic $k$-connected subgraph on more than $(n - k + 1)/2$ vertices for $i = 1,2,3$, where $G_i$ is as defined above. Using this assumption, the proof goes through exactly as above, except inequality (\ref{31kd}) becomes
\begin{equation}\label{31kd'}
0 \; \le \; \sum_i b_i \; \le \; \frac{n - 3k + 3}{2} - 2c \; \le \; \sum_{i,j} x_{ij} - \sum_i z_i - 3k + 3.
\end{equation}
If $x_{ij} + x_{i\ell} - z_i \ge k$ for any triple $\{i,j,\ell\}$, then we would be done as in the proof of Theorem~\ref{31k}, so assume not. So inequality (\ref{31kd'}) is in fact an equality. But then $$\frac{n - 3k + 3}{2} \; = \; 2c$$ with $c \in \N$, which means that $n + k \equiv 1 \pmod 4$, a contradiction. This proves the corollary.
\end{proof}

\begin{rmk}
The bound $n \ge 480k$ is, of course, likely to be far from best possible. By Lemma~\ref{r1kupper}, we know that $n = 6k - 6$ is not sufficient to guarantee the existence of a monochromatic $k$--connected subgraph. We conjecture, along the lines of Bollob\'as and Gy\'arf\'as, that if $n \ge 6k - 5$ then $m(n,3,1,k) \ge (n-k+1)/2$ (this is Conjecture~\ref{conjr1k} in the case $r = 3$).
\end{rmk}

We finish by stating the obvious question: what happens when $r - 1$ is not a prime power? Our lower bound still holds in this case, so we have the following easy corollary of (the proof of) Theorem~\ref{r1k}.

\begin{cor}
Let $r,k \in \N$, and $n \to \infty$. Let $r'$ be the largest integer less than or equal to $r$ such that $r' - 1$ is a prime power. Then
$$\frac{n}{r-1} + o(n) \le m(n,r,1,k) \le \frac{n}{r'-1} + o(n).$$
In particular, $$\left(\frac{1}{6} + o(1)\right)n \le m(n,7,1,k) \le \left(\frac{1}{5} + o(1)\right)n.$$
\end{cor}

\begin{prob}
Find a constant $c = c(r)$ (if one exists) such that $$m(n,r,1,k) = (c + o(1))n$$ for those $r \in \N$ which are not prime powers.
\end{prob}

\section{Acknowledgements}

The authors would like to thank B\'ela Bollob\'as for suggesting the problem to them, and for his ideas and encouragement. They would also like to thank ETH Z\"urich, and Trinity College, Cambridge, where part of this research was carried out.


\begin{thebibliography}{99}

\bibitem{MGT} B.~Bollob{\'a}s, Modern Graph Theory, Springer--Verlag, New York, 1998.

\bibitem{BG} B.~Bollob{\'a}s and A.~Gy{\'a}rf{\'a}s, Highly connected monochromatic subgraphs (manuscript).

\bibitem{ZF} Z. F{\"u}redi, Maximum degree and fractional matchings in uniform hypergraphs, \emph{Combinatorica}, \textbf{1} (1981), 155--162.

\bibitem{AG} A.~Gy{\'a}rf{\'a}s, Partition coverings and blocking sets of hypergraphs (in Hungarian), \emph{Comm. Comp. Automat. Inst. Hung. Acad. Sci.}, \textbf{71} (1977), 62pp.

\bibitem{213} H. Liu, R. Morris and N. Prince, Highly connected monochromatic subgraphs of multicoloured graphs: addendum (manuscript).

\bibitem{HNR2} H. Liu, R. Morris and N. Prince, Highly connected multicoloured subgraphs of multicoloured graphs, submitted.

\bibitem{Mader} W. Mader, Existenz $n$--fach zusammenh{\"a}ngender Teilgraphen in Graphen gen{\"u}gend grosser Kantendichte, \emph{Abh. Math. Sem. Univ. Hamburg}, \textbf{37} (1972), 86--97.

\end{thebibliography}
\end{document}